\documentclass[11pt]{amsart}

\usepackage{epigamath}

\usepackage[english]{babel}

\numberwithin{equation}{section}

\usepackage{tikz-cd}
\usepackage{mathrsfs, mathtools, amssymb, tikz,bm, amscd, scrextend,  verbatim, amsmath, eucal}
\usepackage{microtype}
\usepackage{verbatim}

\usetikzlibrary{matrix, arrows}

\input xy
\xyoption{all} 

\usepackage{esvect}

\newtheorem{theorem}{Theorem}[section]
\newtheorem{lemma}[theorem]{Lemma}
\newtheorem{prop}[theorem]{Proposition}
\newtheorem{proposition}[theorem]{Proposition}

\newtheorem{conditions}{Combinatorial conditions}
\newenvironment{cond}[1]{%
  \renewcommand\theconditions{\rm{#1}}%
  \conditions
}{\endconditions}

\theoremstyle{definition}
\newtheorem{notat}[theorem]{Notation}
\newtheorem{definition}[theorem]{Definition}
\newtheorem{defnote}[theorem]{Definition-Notation}

\theoremstyle{remark}

\newtheorem{question}{Question}
\newtheorem{example}[theorem]{Example}
\newtheorem{remark}[theorem]{Remark}

\newcommand{\bP}{{\mathbb P}}

\newcommand{\bLl}{\mathbb L^{\mathrm{log}}}
\newcommand{\Ext}{\mathrm{Ext}}
\DeclareMathOperator{\Hom}{Hom}

\DeclareMathOperator{\Def}{Def}

\newcommand{\calK}{\mathcal K}
\newcommand{\calO}{\mathcal O}
\newcommand{\calM}{\mathcal M}
\newcommand{\calH}{\mathcal H}
\newcommand{\calF}{\mathcal F}
\newcommand{\calG}{\mathcal G}

\newcommand{\calN}{\mathcal N}

\newcommand{\calX}{\mathcal X}
\newcommand{\calC}{\mathcal C}
\newcommand{\calP}{\mathcal P}

\newcommand{\Supp}{\mathrm{Supp}}
\newcommand{\Spec}{\mathrm{Spec}}
\newcommand{\Spf}{\mathrm{Spf}}

\newcommand{\Pic}{\mathrm{Pic}}

\newcommand{\maps}{\overline{\mathcal{M}}_{0,n}(\bP^1, n)}
\newcommand{\mapso}{\mathrm{Maps}^\circ_n}
\newcommand{\mapsc}{\overline{\mathrm{Maps}}_n}
\newcommand{\mapsd}{\mathcal{M}_{(d_1, \ldots, d_n)}}
\newcommand{\mapsdc}{\overline{\mathcal{M}}_{(d_1, \ldots, d_n)}}

\usepackage[shortlabels]{enumitem}
\setlist[enumerate,1]{label={\rm(\arabic*)}, ref={\rm\arabic*}}

\newcommand{\supth}[1]{\ensuremath{#1^{\mathrm{th}}}}
\renewcommand{\log}{\mathrm{log}}
\newcommand{\qis}{\mathrm{qis}}
\newcommand{\trop}{\mathrm{trop}}
\newcommand{\sm}{\mathrm{sm}}

\EpigaVolumeYear{7}{2023} \EpigaArticleNr{2} \ReceivedOn{November 11, 2021}
\InFinalFormOn{July 23, 2022}
\AcceptedOn{September 17, 2022}

\title{Smoothability of relative stable maps to stacky curves}
\titlemark{Smoothability of relative stable maps to stacky curves}

\author{Kenneth Ascher}
\address{Department of Mathematics, University of California, 
Irvine, CA 92697-3875, USA}
\email{kascher@uci.edu}

\author{Dori Bejleri}
\address{Department of Mathematics, 
Harvard University,
Cambridge, MA 02138, USA}
\email{bejleri@math.university.edu }

\authormark{K.~Ascher and D.~Bejleri}

\AbstractInEnglish{Using log geometry, we study smoothability of genus $0$ twisted stable maps to stacky curves relative to a collection of marked points. One application is to smoothing semi-log canonical fibered surfaces with marked singular fibers.}



\acknowledgement{K.A. is supported in part by funds from NSF Grant DMS-2140781 (formerly DMS-2001408). D.B. is supported in part by an NSF postdoctoral fellowship.}

\begin{document}


\maketitle

\begin{prelims}

\DisplayAbstractInEnglish

\bigskip

\DisplayKeyWords

\medskip

\DisplayMSCclass

\end{prelims}


\newpage

\setcounter{tocdepth}{1}

\tableofcontents


\section{Introduction}

Consider the moduli space $\mapso$ of degree $n$ maps $f\colon C \to \mathbb{P}^1$ 
from a smooth rational curve $C$ such that $f$ is unramified over infinity. Marking the preimage $f^{-1}(\infty) = \{p_1, \ldots, p_n\}$ of such a map induces a locally closed embedding 
\[
\widetilde{\mapso} \hookrightarrow \maps
\]
of the $\mathfrak{S}_n$-torsor $\widetilde{\mapso} \to \mapso$ parametrizing such $f \colon C \to \mathbb{P}^1$ as well as a labeling of $f^{-1}(\infty)$ into the space of $n$-pointed degree $n$ genus $0$ stable maps to $\mathbb{P}^1$. The image of this embedding is the locus of stable maps $(f \colon C \to \mathbb{P}^1, p_1, \ldots, p_n)$ such that $C$ is smooth and $f(p_i) = \infty$ for all $i$. Let $\mapsc$ denote the closure of this locus. 

\begin{question}\label{q1} Is there a combinatorial description of the boundary $\mapsc \setminus \widetilde{\mapso}$? 
\end{question}

There are several natural combinatorial conditions (see Proposition~\ref{prop:necessary} and the discussion preceding it) that are necessary for a stable map $(f \colon C \xrightarrow{n:1} \mathbb{P}^1, p_1, \ldots, p_n)$ to lie inside $\mapsc$:

\begin{enumerate} 
\item The evaluations
  satisfy $f(p_i) = \infty$ for all $i$. 
    \item Every point of $f^{-1}(\infty) \subset C$ either is a marked point or is on a contracted component. 
    \item For each maximal connected closed subvariety $T \subset C$ contracted to $\infty$ by $f$, we have
    \[ \# \textrm{ marked points on $T$ } = \sum \textrm{ ramification index of $f$ at the nodes $T \cap (\overline{C\setminus T})$.}\]
\end{enumerate} 

More generally, for any tuple of non-negative integers $(d_1, \ldots, d_n)$ with $d = \sum d_i$, one can consider the locally closed subset $\mapsd \subset \overline{\mathcal{M}}_{0,n}(\bP^1, d)$ of $n$-pointed degree $d$ genus $0$ maps $(f \colon C \to \bP^1, p_1, \ldots, p_n)$ such that $f(p_i) = \infty$ and $f$ is ramified of order $d_i$ at $p_i$ if $d_i > 0$. If $d_i = 0$, there is no condition imposed on the $f$ at $p_i$. Denoting by $\mapsdc$ the closure of $\mapsd$, we can ask the variant of Question~\ref{q1} for this space. The necessary conditions above naturally generalize to Conditions~\ref{conditions} below. 

\begin{figure}
    \centering
    \includegraphics[scale=.5]{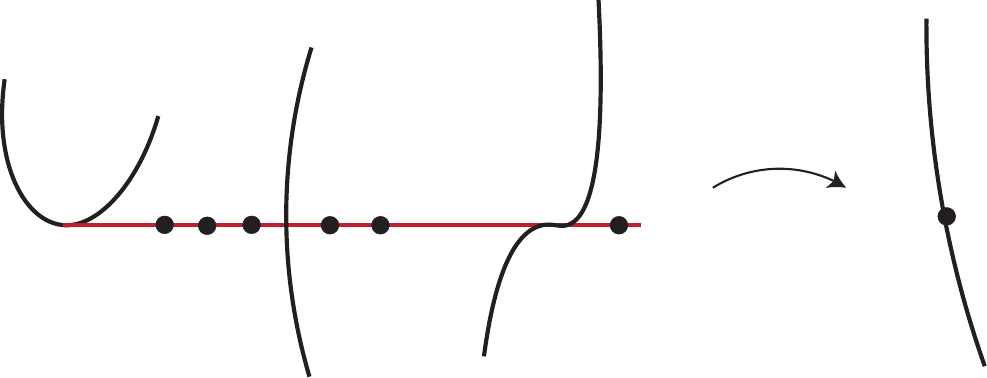}
    \caption{\small{In this example, the contracted component is attached to a double branched, unbranched, and triply branched point, respectively, and thus must have $2 + 1+ 3 = 6$ marked points.}}
    \label{fig:condition3}
\end{figure}

\begin{definition} Let $(f \colon C \to X, p_1, \ldots, p_n)$ be a prestable map to a smooth curve, and fix a tuple of positive integers $(d_1, \ldots, d_n)$ and a point $x \in X$. We say that $f$ is a \textbf{relative map to $(X,x)$ with tangency $(d_1,\ldots, d_n)$} if it satisfies the following conditions:

\begin{cond}{($\ast$)} $($\cite[Definition 1.1 \& Remark 1.7]{Gathmann}$)$\label{conditions} \leavevmode
\begin{enumerate}
    \item\label{cond1} The evaluations satisfy $f(p_i) = x$ for all $i$.
    \item\label{cond2} Every point of $f^{-1}(x) \subset C$ either is a marked point or is on a contracted component. 
    \item\label{cond3} For each maximal connected closed subvariety $T \subset f^{-1}(x)$, we have
    \[ \sum_{p_i \in T} d_i = \sum_{q \in T \cap (\overline{C \setminus T})} e_f(q),\]
    where the first sum is over marked points contained in $T$ and $e_f(q)$ denotes the ramification of $f|_{\overline{C \setminus T}}$ at the point $q$ $($see Figure~\ref{fig:condition3}\,$)$. 
\end{enumerate}
\end{cond} 
\end{definition}

This question in much greater generality was studied by Gathmann~\cite[Proposition 1.14]{Gathmann}, building off previous work of Vakil~\cite[Theorem 6.1]{Vakil}. In particular, Gathmann showed that Conditions~\ref{conditions} relative to $\infty$ are both necessary and sufficient for a stable map to lie in $\mapsdc$. As a consequence, the set of points of $\overline{\mathcal{M}}_{0,n}(\bP^1, d)$ satisfying Conditions~\ref{conditions} relative to $\infty$ are the points of an irreducible closed substack. See also the balancing condition of Gross--Siebert~\cite[Definition 1.12 and Lemma 1.15]{GS}.

In~\cite{av}, Abramovich and Vistoli introduced moduli spaces of \emph{twisted stable maps}, which allows the target to instead be a Deligne--Mumford stack. To form a compact moduli space, the source curves obtain a stacky structure. Let $\calK_{0,n}(\calX, d)$ denote the moduli space of $n$-pointed genus $0$ degree $d$ twisted stable maps to a Deligne--Mumford stack $\calX$. The goal of this paper is to study the analogue of Question~\ref{q1} for genus $0$ twisted stable maps to a weighted projective line $\calP(a,b)$, or more generally a genus $0$ Deligne--Mumford curve.

When $\calX$ is $\calP(a,b)$ and $\infty \in \calP(a,b)$ is a fixed point away from $[0:1]$ and $[1:0]$, we have the following (see Theorem~\ref{mainthm} for a more general result). For this example, the reader can keep in mind   $\overline{\calM}_{1,1} = \calP(4,6)$ and the point $j = \infty$ parametrizing a nodal elliptic curve. 

\begin{theorem}\label{thm:main1} Let $(f \colon \calC \to \calP(a,b), p_1, \ldots, p_n)$ be an $n$-pointed genus $0$ twisted stable map such that the coarse map $(g \colon C \to \mathbb{P}^1, q_1, \ldots, q_n)$ satisfies Conditions~\ref{conditions} with respect to $\infty$. Then $f$ is smoothable in a family with generic fiber satisfying: $f(p_i) = \infty$ and $f$ is ramified to order $d_i$ at $p_i$ for all $i = 1, \ldots, n$. 
\end{theorem}
\noindent In what follows, we set up notation needed to state our more general results.

\medskip
 Let $(\calX, x_1, \ldots, x_r)$ be a smooth and proper $1$-dimensional genus $0$ Deligne--Mumford curve, and suppose the $x_i \in \calX$ are points where the coarse moduli space map is \'etale (equivalently, the $x_i$ have the same stabilizer as the generic point of $\calX$\,). For each $j = 1, \ldots, r$, let $\Gamma_j = (d_{j1}, \ldots, d_{jn_j})$ be a tuple of positive integers, and fix $n_0 \geq 0$. Set $n = \sum_{j=0}^r n_j$ and 
$$
d = \sum_{j=1}^r \sum_{k = 1}^{n_j} d_{jk},
$$
and let $\Gamma = (n_0,\{\Gamma_1,x_1\}, \ldots, \{\Gamma_r,x_r\})$ be the tuple of combinatorial data. 

\begin{defnote}\label{def:MGamma} For any tuple of combinatorial data $\Gamma$, let $\boldsymbol{\mathcal{M}_\Gamma(\calX)}$ be the locally closed substack of $\calK_{0,n}(\calX, d)$ parametrizing $n$-pointed genus $0$ degree $d$ twisted stable maps $(f \colon \calC \to \calX, \{\{p_{jk}\}_{k = 1}^{n_j}\}_{j = 0}^r)$ such that
\begin{enumerate}
    \item\label{def:MGamma-1} the $p_{jk}$ are marked points with stabilizer of order $a_{jk}$ lying over smooth points $q_{jk}$ of the coarse space~$C$,
    \item\label{def:MGamma-2} for each $j > 0$, the image of $p_{jk}$ is $f(p_{jk}) = x_j$,
    \item\label{def:MGamma-3} for each $j > 0$, the coarse map $h \colon C \to X$ is ramified to order $d_{jk}$ at $q_{jk}$, and
    \item\label{def:MGamma-4} $\calC$ is smooth. 
\end{enumerate}

\noindent Similarly, we let $\boldsymbol{\calN_{\Gamma}(\calX)}$ be the locally closed substack of $\calK_{0,n}(\calX,d)$ of maps satisfying conditions~\eqref{def:MGamma-1}, \eqref{def:MGamma-2} and~\eqref{def:MGamma-3} above as well as the condition 
\begin{enumerate}[label={\rm(4${}^\prime$)}, ref=4${}^\prime$]
\item\label{def:MGamma-4'} $\calC$ is smooth in a neighborhood of $f^{-1}(x_j)$ for all $j$.
\end{enumerate}

\noindent We denote by $\boldsymbol{\overline{\calM}_{\Gamma}(\calX)}$ (resp.\ $\boldsymbol{\overline{\calN}_{\Gamma}(\calX)}$) the closure of $\calM_{\Gamma}(\calX)$ (resp.\ $\calN_{\Gamma}(\calX)$) inside $\calK_{0,n}(\calX, d)$. 
\end{defnote} 

Note that here and throughout, we are fixing a bijection between the set $\{1, \ldots, n\}$ and the set $\{ \{(j,k)\}_{k = 1}^{n_j} \}_{j = 0}^r$ indexing the marked points. Note also that the marked points with $j = 0$ are the ones with no tangency conditions, so in the situation of $\mapsd$, we may assume for notational convenience that the first $n_0$ entries $(d_1, \ldots, d_{n_0})$ of $(d_1, \ldots, d_n)$ are $0$ so that $\mapsd = \calM_\Gamma(\mathbb{P}^1)$, where $\Gamma = (n_0, \{(d_{n_0 + 1}, \ldots, d_{n}), \infty\})$.

\begin{defnote}
For any tuple of combinatorial data $\Gamma$, let $\boldsymbol{\calK_\Gamma(\calX)}$ denote the subset of $\calK_{0,n}(\calX, d)$ parametrizing those twisted stable maps $(f \colon \calC \to \calX, \{\{p_{jk}\}_{k = 1}^{n_j}\}_{j = 0}^r)$ such that for each $j = 1, \ldots, r$, the coarse moduli map $(h \colon C \to X, q_{j1}, \ldots, q_{jn_j})$ is a relative map to $(X, x_j)$ with tangency $\Gamma_j$. 
\end{defnote} 

\begin{theorem}\label{mainthm} Let $\calX$ and $\Gamma$ be as above. Then we have an equality $\overline{\calN}_\Gamma(\calX) = \calK_{\Gamma}(\calX)$. That is, every twisted stable map whose coarse moduli map satisfies the relative condition for $\{\Gamma_j, x_j\}$ for each $j = 1, \dots, r$ is smoothable to a twisted stable map parametrized by $\calN_\Gamma(\calX)$. In particular, $\calK_{\Gamma}(\calX)$ is the set of points of a closed substack of~$\calK_{0,n}(\calX, d)$. 
\end{theorem}

\begin{remark} Note that Theorem~\ref{mainthm} is local on the target. Indeed,  the definitions of both $\calK_{\Gamma}(\calX)$ and $\calN_\Gamma(\calX)$ are in terms of local conditions around the points $x_i \in \calX$. Therefore, the theorem has a natural generalization to higher-genus maps to a higher-genus target $\calX$ provided we only consider those maps for which the $1$-dimensional components of the preimages $f^{-1}(x_i)$ are rational curves. In this case, the theorem reads that any such map which is also contained in $\calK_\Gamma(\calX)$ is smoothable in a family with generic fiber contained in~$\calN_\Gamma(\calX)$. 
\end{remark} 

When the target is a weighted projective line, we obtain the following stronger statement. 

\begin{theorem}\label{mainthm2} Let $\calX = \calP(a,b)$ be a weighted projective line, and let $\Gamma$ be a tuple of combinatorial data as above. Then we have an equality $\overline{\calM}_\Gamma(\calX) = \calK_\Gamma(\calX)$. That is, every twisted stable map which satisfies the relative condition for $\{\Gamma_j, x_j\}$ for each $j = 1, \ldots, r$ is smoothable to a family of stable maps from a smooth rational curve satisfying $f(p_{jk}) = x_j$ and with ramification $d_{jk}$ at $p_{jk}$. 
\end{theorem}

\subsection{Applications to moduli of fibered surfaces} 

Our original motivation for writing this paper came from studying compactifications of the moduli space of fibered surfaces. Twisted stable maps are used in~\cite{avfibered, tsm} to construct a compactification $\calF_{g,n}^v(\gamma, \nu)$ of the moduli space of genus $\gamma$ fibrations over a genus $g$ curve with $\nu$ marked sections and $n$ marked fibers. The objects of the boundary are certain semi-log canonical unions of birationally fibered surfaces called \emph{twisted surfaces}. This compactification is closely related to the compactification via stable log varieties from the minimal model program. In~\cite[Section 1.4]{tsm}, we proposed the problem of using log geometry to give the main component a moduli-theoretic interpretation and classify the boundary components. The present paper solves this problem for elliptic fibrations with marked singular fibers (we refer the reader to~\cite{master} and~\cite[Section 4]{k3}). A key observation is that Conditions~\ref{conditions} relative to $\infty \in \overline{\calM}_{1,1}$ as well as the choice of stabilizers on the marked points translate to conditions on the configuration of singular fibers on the components of the twisted elliptic surface; see \cite[Propositions 4.1 and 4.4]{k3}. 

\begin{theorem} Theorem~\ref{mainthm2} gives necessary and sufficient combinatorial conditions for a twisted elliptic surface over a genus $0$ curve with marked singular fibers to be smoothable to an elliptic surface over $\mathbb{P}^1$ with marked singular fibers.
\end{theorem}

For convenience, we work over an algebraically closed field of characteristic $0$. 

\subsection*{Acknowledgements} We thank Dan Abramovich, Qile Chen, Martin Olsson, Dhruv Ranganathan, Ravi Vakil, Jonathan Wise, and David Yang for helpful discussions.  We thank the referees for their helpful comments.

\section{Genus 0 relative stable maps to \texorpdfstring{$\boldsymbol{(\mathbb{P}^1, \infty)}$}{(P\textasciicircum 1, infinity)}}\label{sec:combinatorics}

In this section, we prove the special case of Theorem~\ref{thm:main1} where the target is $\bP^1$; \textit{i.e.}, $a = b = 1$. This special case was originally proved by Gathmann \cite{Gathmann}. 
Our approach differs from that of~\cite{Gathmann} in that we give a direct construction of smoothings of comb-type maps (see Propositions~\ref{prop:sufficient} and~\ref{prop:linequiv1}) rather than appealing to~\cite[Theorem 6.1]{Vakil}. This will be a key step in the proof of the general case of Theorem~\ref{thm:main1}.  For the remainder of this section, fix positive integers $\Gamma = (d_1, \ldots, d_n)$ with $d = \sum d_i$.

\begin{figure}
    \centering
    \includegraphics[scale=.5]{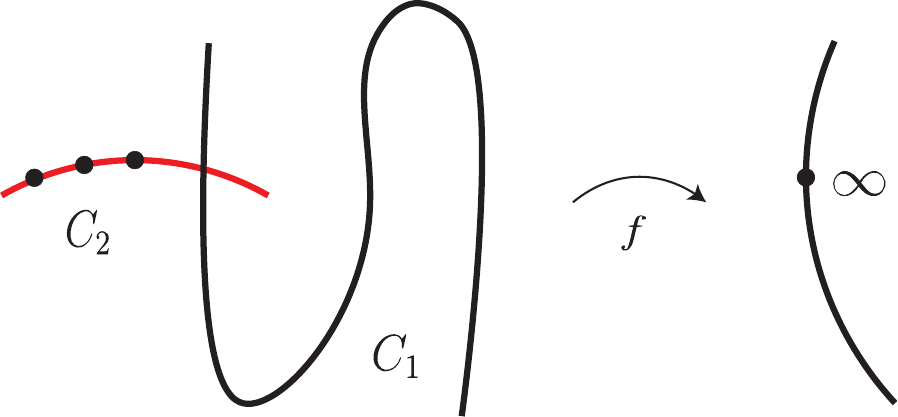}
    \caption{\small{Maps of this type are given by choosing the map $f_1 = f|_{C_1}$ (which is parametrized by $\mapsc$), a point of the finite set $f_1^{-1}(\infty)$ along which to attach $C_2$, and the configuration of the $n + 1$ special points on $C_2$ (which is parametrized by $\mathcal{M}_{0,n+1}$), yielding the dimension count $\dim \mapsc + n - 2$.}}
    \label{fig:example1}
\end{figure}

\begin{example}\label{example1} 
We begin with some examples motivating Conditions~\ref{conditions} when $\Gamma = (1, \ldots, 1)$. First, it is clear that Condition~\ref{conditions}\eqref{cond1} is required as the evaluation condition is closed. However, consider a degree $n$ map $f_1\colon (C_1,q_1) \to \bP^1$ from a smooth rational curve $C_1$ with $f(q_1) = \infty$. Let $C = C_1 \cup_{q_1, q_2} C_2$ be a nodal union of two rational curves, and let $p_1, \ldots, p_n \in C_2 \setminus q_2$ be $n$ marked points. Then there is an $n$-marked degree $n$ stable map $(f \colon C \to \bP^1, p_1, \ldots, p_n)$ given by taking the map $f_1$ on $C_1$ and the constant map with image $\infty$ on $C_2$. While this map satisfies Condition~\ref{conditions}\eqref{cond1}, a simple dimension count shows that the dimension of this locus inside $\maps$ is equal to $\dim \mapsc + n -2$, and so this condition alone is not enough to cut out the locus $\mapsc$ at least for $n \geq 3$ (see Figure~\ref{fig:example1}).
\end{example}

\begin{figure}
    \centering
    \includegraphics[scale=.5]{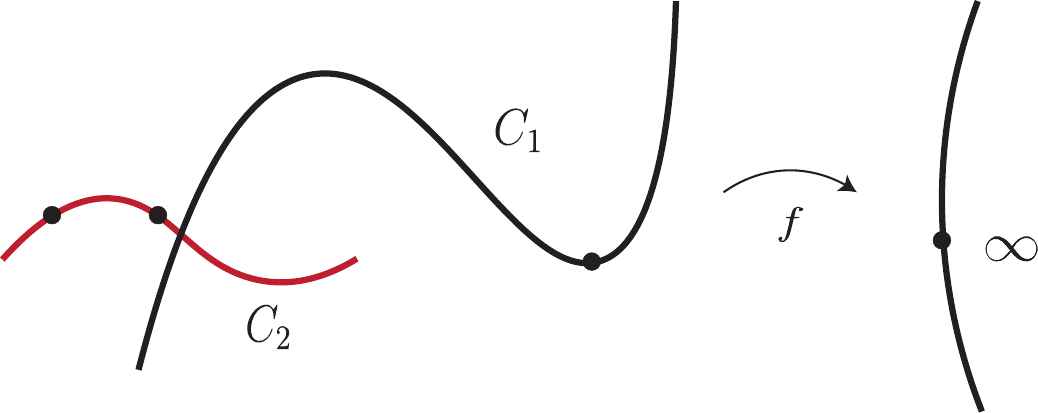}
    \caption{\small{The map in Example~\ref{example2} which satisfies Conditions~\eqref{cond1} and~\eqref{cond2} but not~\eqref{cond3} since the map is unramified where the contracted component $C_2$ is attached and ramified over $\infty$ elsewhere.}}
    \label{fig:example2}
\end{figure}

The above example motivates Condition~\ref{conditions}\eqref{cond2}, which in this case requires that the map $f_1$ be totally ramified at $q_1$ so that there are no other points of $C_1$ in $f^{-1}(\infty)$. Requiring $f_1$ to be totally ramified at $q_1$ means that we impose the vanishing of $n$ derivatives, which is a codimension $n$ condition on $\mapsc$. Thus the locus of maps of this combinatorial type satisfying both Conditions~\ref{conditions}\eqref{cond1} and~\ref{conditions}\eqref{cond2} has dimension $\dim \mapsc - 2$. Note that such maps automatically satisfy Condition~\ref{conditions}\eqref{cond3} as well. However, the following example illustrates that Conditions~\ref{conditions}\eqref{cond1} and~\ref{conditions}\eqref{cond2} do not imply Condition~\ref{conditions}\eqref{cond3} in general.

\begin{example}\label{example2}
Consider a degree $3$ map $f_1\colon (C_1,q_1, q_2) \to \bP^1$ from a smooth rational curve such that $f(q_i) = \infty$, and suppose that $f$ has ramification index $i$ at $q_i$. Let $(C_2, q_1', p_1, p_2)$ be a smooth pointed rational curve, and let $C = C_1 \cup_{q_1, q_1'} C_2$. Then we have a $3$-pointed degree $3$ stable map $(f \colon C \to \mathbb{P}^1, q_2, p_1, p_2)$ where $f|_{C_1} = f_1$ and $f|_{C_2}$ is constant $\infty$. This map satisfies Conditions~\ref{conditions}\eqref{cond1} and~\ref{conditions}\eqref{cond2} but not~\ref{conditions}\eqref{cond3}. Note furthermore that this stable map cannot lie in $\mapsc$ as $f_1$ is ramified to order $2$ at the marked point $q_2$, but the only way this can happen in the limit of a family in $\mapso$ is if two marked points collided at $q_2$.  \end{example}

\begin{remark} Note that a degenerate special case of Condition~\ref{conditions}\eqref{cond3} is when $T = q$ is itself a point. Since~$T$ must be a maximal closed subvariety contracted by $f$, this means that $q$ cannot lie on a contracted component. Therefore, by Condition~\ref{conditions}\eqref{cond2}, $q$ must be a marked point. Then Condition~\ref{conditions}\eqref{cond3} reads that $f$ must be unramified at $q$. Indeed, we saw this was necessary in the above Example~\ref{example2} (see Figure~\ref{fig:example2}). 

\end{remark}

We now prove that conditions~\eqref{cond1},~\eqref{cond2}, and~\eqref{cond3} of Conditions~\ref{conditions} are necessary. 

\begin{prop}\label{prop:necessary}
Conditions~\ref{conditions}~\eqref{cond1},~\eqref{cond2}, and~\eqref{cond3} are necessary for a stable map to lie in $\mapsdc$. 
\end{prop}

\begin{proof}
It is clear that~\eqref{cond1} is necessary as evaluation is continuous so the condition $f(p_i) = \infty$ is closed, so we proceed to~\eqref{cond2} and~\eqref{cond3}. Consider a $1$-parameter family of stable maps 
\[ \begin{tikzcd}
C \arrow{r}{f}  \arrow[d, "\pi"] & \bP^1 \\ 
  S \arrow[u, bend left, "\sigma_i"] 
\end{tikzcd} \] 
over the spectrum  $S$ of a DVR  with generic fiber lying in $\mapso$, and denote by $(f_0 \colon C_0 \to \bP^1, p_1, \ldots, p_n)$ the central fiber. 

First we show~\eqref{cond2}. Let $U \subset C$ denote the open complement of the locus of components contracted by the map $f$, and consider $D := f^{-1}(\infty)$. As the total space $C$ is normal and $f$ is non-constant, $D$ is a Cartier divisor. Moreover, note that any $\pi$-vertical component of $D$ is contracted by $f$ and therefore is not in $U$. Thus, restricting to $U$, we see that $D|_U$ is a Cartier divisor which is horizontal over $S$, and therefore any point of $D|_U$ lying over $0 \in S$ must be in the closure of the marked points of the generic fiber.

Finally, we show~\eqref{cond3}. Ramification corresponds to a polynomial having a multiple root at a point, and this multiplicity takes into account precisely how many points collided, \textit{i.e.}, the number of marked points on the relevant contracted component. More formally, we can consider the intersection product $T\cdot D$. By the projection formula, $T\cdot D = 0$. On the other hand, $D = \sum \sigma_i + E$, where $E$ is the $\pi$-vertical component. Computing $T\cdot (\sum \sigma_i + E) = 0$ in terms of local multiplicities gives exactly the equality in Condition~\ref{conditions}\eqref{cond3}. 
\end{proof}

Our task now is to show that Conditions~\ref{conditions}\eqref{cond1}--\eqref{cond3} are sufficient for a stable map $(f\colon  C \to \bP^1, p_1, \ldots, p_n)$ to lie in $\mapsdc$. We will do this by constructing a smoothing of the marked curve and a linear series on the total space which restricts to $f$ on the central fiber but whose generic fiber lies in $\mapsd$.  

We begin with a preliminary lemma that shows that nodal curves can be smoothed to surfaces admitting $A_m$ singularities for any $m$ at the nodes.  

\begin{lemma}\label{lemma:ansing}
Let $C_0$ be a genus $0$ nodal curve with $k$ nodes $q_1, \ldots, q_k$, and let $m_1, \ldots, m_k \geq 1$. Then there exists a smooth $C \to S$, where $S = \Spec(R)$ is the spectrum of a DVR, such that the total space $C$ has an $A_{m_i-1}$ singularity at $q_i$ for each $i = 1, \ldots, k$. 
\end{lemma}

\begin{proof}
The miniversal deformation space of a pointed prestable curve of genus $0$ is $k \llbracket t_1, \dots, t_k \rrbracket$, where $t_i$ is the smoothing parameter of the $\supth{i}$ node; \textit{i.e.}, formally locally around the $\supth{i}$ node, the miniversal family of curves looks like $xy=t_i$ inside $\mathbb{A}^2_{k \llbracket t_1, \dots, t_k \rrbracket}$. Let $R = \llbracket z \rrbracket$ be the DVR in the statement. To prove the lemma, it suffices to note we can construct a map of rings  $k \llbracket t_1 \dots t_k \rrbracket \to k \llbracket z \rrbracket$ such that $t_i \mapsto z^{m_i}$. Then formally locally around the $\supth{i}$ node, the family of curves over $\Spec(R)$ will be isomorphic to $xy = z^{m_i}$, as required. 
\end{proof} 

We now show that it suffices to consider the case of a stable map where every maximal connected subtree contracted by $f$ (as in Condition~\ref{conditions}\eqref{cond3}) is irreducible. 

\begin{lemma}\label{lem:comb}
Let $(f_0\colon C_0 \to \bP^1, p_1, \ldots, p_n)$ be a stable map satisfying Conditions~\ref{conditions}\eqref{cond1},~\eqref{cond2}, and~\eqref{cond3}. Then there exists a deformation to a family of stable maps $(f\colon C \to \bP^1, \sigma_i) \to \Spec(R)$ over the spectrum of a DVR such that 

\begin{enumerate}
\item the generic fiber $(f_\eta, C_\eta \to \bP^1, (\sigma_i)_\eta)$ satisfies Conditions~\ref{conditions}\eqref{cond1},~\eqref{cond2}, and~\eqref{cond3}, and
\item every connected component of $f_\eta^{-1}(\infty) \subset C_\eta$ is irreducible. 
\end{enumerate}
\end{lemma}

\begin{proof} Let $E^1, \ldots, E^k$ be the $1$-dimensional connected components of $f_0^{-1}(\infty)$, and write $C^1, \ldots, C^r$ for the connected components of the closure of the complement $C_0 \setminus \{E^i\}_{i=1}^k$. Each $E^i$ is pointed by $(q_{ij}, p_{il})$, where the $p_{il}$ are the marked points that lie on $E^i$ and the $q_{ij}$ are the points of $E^i$ along which $E^i$ is glued to the non-contracted components of $C$. Then each $E_i$ is a pointed genus $0$ prestable curve, and there exists a smoothing $\widetilde{E}^i \to \Spec(R)$ with sections $\tau_{ij}$ smoothing the $q_{ij}$ and $\sigma_{il}$ smoothing the $p_{il}$. Now consider the constant family $\bigsqcup_l C^l_R$ with constant sections corresponding to the marked points and the points $q'_{ij}$ along which the $q_{ij}$ are glued. 

We can glue $\bigsqcup_l C^l_R$ with $\bigsqcup_t \widetilde{E}^t$ by identifying the sections $q'_{ij} \times \Spec(R)$ with $\tau_{ij}$ for all $(i,j)$ and call the result $C \to \Spec(R)$, which is now a family of pointed curves. By construction, this is a partial smoothing of the pointed curve $(C_0, p_1, \ldots, p_n)$. Moreover, there is a stable map $f\colon  C \to \bP^1$ constructed by taking $f_R|_{C^l_R} \colon C^l_R \to \bP^1$ for each $l$ and taking the constant map $\infty$ on each $\widetilde{E}^t$. This descends to a map $f$ as desired since $f(q'_{ij}) = \infty$. Conditions~\ref{conditions}\eqref{cond1}, \eqref{cond2} are satisfied by $f_\eta$ by construction, and~\ref{conditions}\eqref{cond3} is satisfied since the maximal irreducible components $\tilde{E}^i_\eta$ of $f_\eta^{-1}(\infty)$ have the same number of marked points as $E^i$ and $f_\eta$ has the same ramification at $(q'_{ij})_\eta$ as $f_0$ at $q'_{ij}$. 
\end{proof}

We now restrict ourselves to the case where the connected components of $f^{-1}(\infty)$ are irreducible. Let $(f \colon C_0 \to \bP^1, p_1, \ldots, p_n)$ be such a map, suppose that $E_1, \ldots, E_k$ are the $1$-dimensional connected components of $f^{-1}(\infty)$, and write $E = \sqcup_i E_i$. By the above reduction, each $E_i$ is a smooth rational curve.

We write the closure of the complement of $E$ in $C_0$ as a union of connected components $C_1, \ldots, C_r$. Thus each $C_i$ is a tree of rational curves such that $(f_0)|_{C_i}$ is non-constant and such that the preimage of $\infty$ is $0$-dimensional. For each $i$, let $I_i \subset \{1, \ldots, n\}$ be the subset of indices $k$ such that $p_k$ lies on $E_i$. If $E_i$ and $C_j$ intersect, we let $q_{ij} \in E_i$ and $t_{ij} \in C_j$ be the points of $E_i$ and $C_j$, respectively, at which they are glued, and we let $e_{ij}$ be the ramification index of $(f_0)|_{C_j}$ at $t_{ij}$. Note that Condition~\ref{conditions}\eqref{cond3} then reads that for each $i$, 
$$
n_i := \sum_{k \in I_i} d_k= \sum_{j} e_{ij}, 
$$
where the left side is a definition and the right sum is over the $j$ such that $E_i$ meets $C_j$. 

Consider a smoothing $(C \to S, \sigma_1, \ldots, \sigma_n)$ of $(C_0, p_1, \ldots, p_n)$ over the spectrum of a DVR, and let $D = \sum d_i \sigma_i$ denote the divisor of $\Gamma$-weighted marked sections on this surface. We wish to construct a rank $1$ linear series on such a smoothing so that the central fiber agrees with $f_0$ and the generic fiber satisfies that $D|_\eta$ is the preimage of $\infty$. Our strategy is to consider the map $\pi \colon C \to C'$ which contracts $E$ to a point and instead construct the appropriate linear series on $C'$. 

\begin{prop}\label{prop:sufficient} In the setting above, there exist a smoothing $(C \to S, \sigma_1, \ldots, \sigma_n)$ and a Cartier divisor $D'$ on~$C'$, where $\pi \colon C \to C'$ is the contraction of the $E_i$
  such that
\begin{enumerate} 
\item $\pi^*D' = D + \sum a_iE_i$ for some $a_i$, where $D$ is the divisor of marked sections, and
\item $a_iE_i\cdot C_j = e_{ij}t_{ij}$ for all $i$ and $j$ such that $E_i$ meets $C_j$. 
\end{enumerate}
\end{prop}

\begin{proof}

For each $i$, we let $a_i = \prod_{j = 1}^r e_{ij}$ and let $m_{ij} = \prod_{k \neq j} e_{ik}$. By Lemma~\ref{lem:comb}, there exists a smoothing $(C \to S, \sigma_1, \ldots, \sigma_n)$ of $(C_0, p_1, \ldots, p_n)$ such that $C$ has an $A_{m_{ij} -1}$ singularity at the node $t_{ij}$. Computing intersection products in the surface $C$, we have $E_i\cdot C_0 = 0$ since $C_0$ is a fiber containing $E_i$. Moreover, if $E_i$ and $C_j$ intersect, then
$$
E_i\cdot C_{j} = \frac{1}{m_{ij}} t_{ij}
$$
since locally around the node $t_i$, the curves $E_i$ and $C_j$ are distinct lines through an $A_{m_{ij} - 1}$ singularity. Therefore, 
$$
E_i^2 = \sum_j -\frac{1}{m_{ij}}, 
$$
where the sum is over those $j$ for which $C_i$ and $E_j$ meet. Here, we go back and forth between viewing these intersection products as numbers or as divisors on the curves depending on whether it is convenient to emphasize the particular intersection points. 

Now consider the divisor $D + \sum a_iE_i$. We can compute that $D\cdot E_i = \sum_{k \in I_i} d_k = n_i$ since $D$ is the divisor of weighted marked sections and $E_i$ contains $p_k$ for $k \in I_i$. On the other hand,
$$
aE_i\cdot C_j = \frac{a_i}{m_{ij}} = \frac{\prod_j e_{ij}}{\prod_{k \neq j}e_{ik}} = e_{ij}.
$$
Finally, 
$$
\Big(D + \sum a_kE_k\Big)\cdot E_i = n_i - \sum_j \frac{a_i}{m_{ij}} = n_i - \sum_j \frac{\prod_i e_{ij}}{\prod_{k \neq j}e_{ik}} = n_i - \sum_j e_{ij} = 0
$$
by Condition~\ref{conditions}\eqref{cond3}. Here we have used that $E_k\cdot E_i = 0$ if $i \neq k$. 

Therefore, we need to show that $D + \sum a_iE_i$ descends to a Cartier divisor $D'$ along the contraction $\pi \colon C \to C'$. Note that $\pi$ exists, $C'$ is a normal quasiprojective surface, and $\pi_*\calO_C = \calO_{C'}$ by~\cite[Theorem~3.7]{km}. Moreover, we have an exact sequence 
\begin{equation}\label{eqn:exact}
\xymatrix{
0 \ar[r] & \Pic(C') \ar[r]^{\pi^*} & \Pic(C) \ar[r] & \mathbb{Z}, 
}
\end{equation}
where the map $\Pic(C) \to \mathbb{Z}$ is the restriction $L \mapsto L|_E$; see \cite[Corollary 3.17]{km}. In particular, since  $\calO_C(D + \sum a_iE_i)$ is in the kernel of this map,  there exists a line bundle $L$ on $C'$ such that $\pi^*L = \calO_C(D + aE)$. On the other hand, $\pi_*\calO_C = \calO_{C'}$, and so by the projection formula, 
$$
\pi_*\calO_C\Big(D + \sum a_iE_i\Big) = L.
$$
Therefore, there exists a section $s \in H^0(L)$ which pulls back to the section cutting out $D + \sum a_iE_i$ on $C$, and so $D' = \mathrm{div}(s)$ is the required Cartier divisor. \end{proof} 

\begin{remark}\label{remark:Am} Fix a list of positive integers $r_{ij}$ for $i$ and $j$ such that $E_i$ meets $C_j$. Then in the construction of of the smoothing in Proposition~\ref{prop:sufficient}, we can replace $e_{ij}$ with $r_{ij}e_{ij}$ in the definitions of $a_i$ and $m_{ij}$. In this way, we get a smoothing satisfying the properties of Proposition~\ref{prop:sufficient}, where the total space has $A_{m_{ij}-1}$ singularities with $m_{ij}$ divisible by any fixed choice of integers. 
\end{remark}

Let $0 \in \bP^1$ be a general point of the target, and consider the divisor $B_0 \subset C_0$ given by $f_0^{-1}(0)$. Then $B_0$ consists of a union of points on $C_0$ which are disjoint from $E$ and contained in the locus where $f_0$ is \'etale. Since $C \to S$ is a family of genus $0$ curves birational to $\bP^1_S$, we can extend the collection of points $B_0$ into a Cartier divisor $B$ on $C$ which is horizontal and satisfies $B|_{C_0} = B_0$. 

\begin{prop}\label{prop:linequiv2}
With notation as above, $D + \sum a_iE_i \sim B$ are linearly equivalent Cartier divisors on $C$. 
\end{prop}

\begin{proof}
  Let $\psi\colon C \to S$ denote the morphism to the base, and consider $\psi_*N$, where $N = \calO_{C}(D + \sum a_iE_i - B)$. First we claim that $H^1(C_0, N|_{C_0}) = 0$. Indeed, since $C_0$ is a nodal genus $0$ curve, it suffices to show that the degree of $N$ restricted to each component of $C_0$ is at least $0$.
  Now $N|_E = \calO_E$ by construction since $(D + \sum a_iE_i)\cdot E = 0$ and $B$ avoids $E$. On the other hand, for each $i$, we have $(D + \sum a_iE_i)\cdot C_j = f_0^{-1}(\infty)|_{C_j}$ and $B\cdot C_j = f_0^{-1}(0)|_{C_j}$, so these divisors are linearly equivalent. Therefore, we in fact have that $N$ is trivial on each component of $C_0$. Hence, by cohomology and base change, $\psi_*N$ is a vector bundle whose formation commutes with base change. Moreover, over the generic fiber, $C_\eta \cong \bP^1_\eta$, and $B$ and $D + \sum a_iE_i$ are divisors of the same degree $n$ by construction. Therefore, $H^0(C_\eta, N|_{C_\eta}) = 1$. Since any line bundle on the spectrum of a DVR is trivial, we conclude that $\psi_*N \cong \calO_S$. Consequently, $N$ has a non-vanishing section which exhibits the required linear equivalence $D+ \sum a_i E_i\sim B$. \end{proof} 

\begin{remark}
Note that the vanishing claimed in the above proof does not hold if we replace $C_0$ with a genus~$0$ Deligne--Mumford stack. See Remark~\ref{rem:needlog} for an example.
\end{remark}

Now let $B'$ be the image of $B$ under $\pi$, which is also a Cartier divisor since $B$ is contained in the locus where $\pi$ is an isomorphism. By the exact sequence~\eqref{eqn:exact}, the linear equivalence in Proposition~\ref{prop:linequiv2} is equivalent to a linear equivalence $D' \sim B'$. In particular, $D'$ and $B'$ form a basepoint-free rank $1$ linear subseries of $H^0(C', L)$, where $L = \calO_{C'}(D')$.

\begin{prop}\label{prop:linequiv1}
With notation as above, let $g\colon  C' \to \bP^1$ be the morphism induced by the basepoint-free linear series $\langle B', D' \rangle$. Then the composition $C \xrightarrow{\pi} C' \xrightarrow{g} \bP^1$ satisfies $(g \circ \pi)|_{C_0} = f_0$. 
\end{prop}

\begin{proof}
The map $(g \circ \pi)$ contracts $E$ to $\infty$ by construction. Moreover, on each $C_j$, the restrictions of $B'$ and~$D'$ satisfy 
$$
B'|_{C_j} = f_0^{-1}(0)|_{C_j}
$$
and 
$$
D'|_{C_j} = f_0^{-1}(\infty)|_{C_j}
$$
by Proposition~\ref{prop:sufficient}. Therefore, $(g \circ \pi)|_{C_j} = (f_0)_{C_j}$ for all $j = 1, \ldots, k$. We conclude that $(g \circ \pi)|_{C_0} = f_0$ as it agrees with $f_0$ on each component of $C_0$. \end{proof}

Putting these together, we conclude the following. 

\begin{theorem}\label{thm:sufficient} Let $(f_0  \colon C_0 \to \bP^1, p_1, \ldots, p_n)$ be a stable map satisfying Conditions~\ref{conditions}\eqref{cond1},~\eqref{cond2}, and~\eqref{cond3}. Then there exist a smoothing $(C \to S, \sigma_1, \ldots, \sigma_n)$ over $S = \Spec(R)$, the spectrum of a DVR, and a stable map $f \colon C \to \bP^1$ such that the generic fiber $(f_\eta, C_\eta \to \bP^1_\eta, (\sigma_i)_\eta)$ is contained in $\mapsd$ and $f|_{C_0} = f_0$. In particular, Conditions~\ref{conditions}\eqref{cond1},~\eqref{cond2}, and~\eqref{cond3} are both necessary and sufficient for a stable map to be contained in $\mapsdc$. 
\end{theorem}

\begin{proof} By Propositions~\ref{prop:sufficient} and~\ref{prop:linequiv1}, there exist a smoothing $(C \to S, \sigma_i)$ and a map $f \colon C \to \bP^1$ such that $f|_{C_0} = f_0$ and $f^{-1}(\infty) = D + \sum a_iE_i$, where $D$ is the divisor of weighted marked sections and the $E_i$ are contained on $C_0$. Therefore, $f|_{C_\eta}$ is in $\mapsd$. 
\end{proof}

\section{Log twisted curves and maps}

We refer the reader to~\cite{HoM} for the standard definitions in log geometry and to~\cite[Section 2]{logGLSM} for log twisted stable maps. We always work with fine and saturated (fs) log structures. We will use $M_X$ to denote the log structure of $X$, which will be implicit if no confusion arises. The characteristic of the log structure will be denoted by $\overline{M}_X$. 

\subsection{Log geometry}\label{sec:log}

We will need the following results regarding log smooth morphisms. 
\begin{remark}\leavevmode
\begin{itemize}
\item If a morphism of log schemes $f\colon  X \to {\rm pt}$ is log smooth (where the point has the trivial log structure), then $U \subset X$ is a toroidal embedding, where $U$ is the locus with trivial log structure and $M_X$ is the divisorial log structure.
    \item If $C$ is a curve, then log smooth is equivalent to nodal; see~\cite[Theorem 5]{HoM}.
    \item Given a family $f \colon  C \to S$ of nodal curves, there is a \emph{minimal} log structure $M_S$ on $S$ such that any other log structure that makes $f$ log smooth is pulled back from it; see~\cite[Section 7.3]{HoM}. We call the corresponding log structure and structure morphism $f^\flat \colon  f^*M_S \to \tilde{M}_C$ the \emph{canonical} log structure of
      $f \colon  C \to S$. Note here that $\tilde{M}_C$ does not include marked points. If $f$ is equipped with marked points, then we denote the natural log structure by 
    $$
    M_C := \tilde{M}_C \bigoplus_{\calO_C^*} \left( \oplus_{\calO_C^*, i}N_i\right), 
    $$
    where the sum is over marked points and $N_i$ is the divisorial log structure associated to the $\supth{i}$ marked point; we call $M_C$ the \emph{canonical log structure} associated to a pointed curve if there is no confusion. 
\end{itemize}
\end{remark}

There exists a log cotangent complex $\bLl_{X/Y}$ (see~\cite{Olssonlogct} or~\cite[Section 7]{HoM}) for morphisms of log schemes $f\colon  X \to Y$, and deformation theory of log schemes is controlled by the log cotangent complex (see~\cite[Theorem 5.2]{Olssonlogct}).

\begin{remark}\label{rmk:logcotangent}
  If $f$ is a log smooth morphism, then the log cotangent complex $\bLl_{X/Y}$ is represented by the sheaf of log differentials (see \cite[Section~1.1(iii)]{Olssonlogct}).
  There does not exist a distinguished triangle in general; however, Olsson constructs a distinguished triangle for log flat or integral morphisms (see~\cite[1.3]{Olssonlogct}).
\end{remark}

\subsection{(log) Twisted curves}
Our main reference is~\cite{Olsson}. To compactify moduli spaces of maps $f\colon  C \to \calM$, where $\calM$ is a Deligne--Mumford stack, one needs to allow $C$ to be a stack as well, known as a \emph{twisted curve} (see \textit{e.g.}~\cite{av, tsm}).

\begin{definition}
A \emph{twisted curve} is a purely 1-dimensional Deligne--Mumford stack $\calC$, with at most nodes as singularities, satisfying the following conditions:
\begin{enumerate}
    \item If $\pi\colon  \calC \to C$ denotes the coarse space morphism, then $\calC^{\sm} = \pi^{-1}C^{\sm}$, and $\pi$ is an isomorphism over a dense open subset of $C$. 
    \item If $\bar{x} \to C$ is a node such that the strictly henselian local ring $\calO_{C, \bar{x}}$  is the strict henselization of $k[x,y]/(xy)$, then
    \[ \calC \times_C \Spec \left(\calO_{C, \bar{x}}\right) \cong \bigg[\Spec \left(\calO_{C, \bar{x}}[z,w]/(zw, z^m-x, w^m-y)\right) / \mu_m\bigg],\]
    where $\xi \in \mu_m$ acts by $(z,w) \mapsto (\xi z, \xi^{-1}w)$.
\end{enumerate}
\end{definition}

\begin{definition}
An \emph{$n$-pointed twisted curve} $\calC/S$ marked by disjoint closed substacks $\{ \Sigma_i\}_{i=1}^n$ in $\calC$ is assumed to satisfy the following:
\begin{enumerate}
    \item The $\Sigma_i$ are contained in $\calC^{\sm}$. 
    \item Each $\Sigma_i$ is a tame \'etale gerbe over $S$.
    \item The map $\calC^{\sm} \setminus \cup \Sigma_i \to C$ is an open embedding.
\end{enumerate}
\end{definition}

Let $\calX$ be a Deligne–Mumford stack. We say that a fine log structure $M_{\calX}$ is \emph{locally free} if for every geometric point $\bar{x} \to \calX$, the characteristic sheaf satisfies $\overline{M}_{\calX} \cong \mathbb{N}^r$ for some $r$. 

\begin{definition}\label{def:simple}
In the above situation, we say that a morphism of sheaves of monoids $M \to M'$ is \emph{simple} if for every geometric point $\bar{x} \to \calX$, we have
\[\begin{tikzcd}
    \overline{M}_\calX \ar{r} \ar{d} & \overline{M}'_\calX \ar{d} \\
    \mathbb{N}^r \ar[r, "\phi"] & \mathbb{N}^r\rlap{,}
\end{tikzcd}\] 
where $\phi$ is given by $(m_1, \dots, m_r)$.
\end{definition}

\begin{definition}
An \emph{$n$-pointed twisted log curve over $S$} is the data
\[ (C/S, \{\sigma_i, a_i\}, l\colon M_S \to M'_S),\]
where 
\begin{itemize}
    \item $(C, \{\sigma_i\})/S$ is an $n$-pointed nodal curve,
    \item $M_S$ is the minimal log structure for the family $C \to S$,
    \item the $a_i\colon  S \to \mathbb{Z}_{>0}$ are locally constant, and
    \item $l$ is a simple morphism.
\end{itemize}
\end{definition}

\begin{theorem}[\textit{cf.} \protect{\cite[Theorem 1.8]{Olsson}}]\label{thm:olsson}
The fibered category of $n$-pointed twisted curves is naturally equivalent to the stack of $n$-pointed log twisted curves. \end{theorem}

There is a natural map from the stack of twisted curves to the stack of (pre)stable curves induced by taking the coarse space $\pi \colon  \calC \to C$. The induced map on local deformation spaces $\Def(\calC) \to \Def(C)$ can be described as a root stack of order $m_i$ along the boundary divisor $\{t_i = 0\}$, where $t_i$ is the deformation parameter of the $\supth{i}$-node of $C$ and $m_i$ is the stabilizer order of $\calC$ at the $\supth{i}$ node. The stabilizer orders $m_i$ correspond via Theorem~\ref{thm:olsson} to the data of the simple extension $l$ as in Definition~\ref{def:simple}. For more details, see the discussion following~\cite[Theorem 1.9 and Remark 1.10]{Olsson}. 

\begin{example} \label{rem:sing} Consider a smoothing of a nodal curve as in Lemma~\ref{lemma:ansing}, and let $(C, M_C) \to (S, M_S)$ be the minimal log structure. The appearance of $A_{m_i-1}$ singularities on the total space of the smoothing at the nodes of $C_0$ is equivalent to the existence of a simple extension $M_S \hookrightarrow M_S'$ with $\phi = (m_1, \ldots, m_k)$. Thus, by Olsson's Theorem~\ref{thm:olsson}, such a smoothing is the coarse space of a smoothing of the twisted curve with stabilizer $\mu_{m_i}$ at the $\supth{i}$ node of $C_0$. Note that we can see this directly by taking the canonical stack of the total space $C$ which introduces a stabilizer  $\mu_{m_i}$  at the $A_{m_i - 1}$ singularity. 
\end{example} 

Let us briefly explicate the construction of Theorem~\ref{thm:olsson} as we will use it further on. If $\pi \colon  \calC \to C$ is the coarse map of a twisted curve over $S$, then the simple extension $l \colon  M_S \hookrightarrow M_S'$ is the map between the minimal log structures of $h \colon  C \to S$ and $f \colon   \calC \to S$, respectively. The pushforward log structure $\pi_*\tilde{M}_\calC$ sits in a pushout square
$$
\xymatrix{h^*M_S \ar[r]^{h^*l} \ar[d]_{h^\flat} & h^*M_{S}' \ar[d]^{\pi_*f^\flat} \\ \tilde{M}_C \ar[r] & \pi_*\tilde{M}_{\calC}\rlap{.} }
$$
Each $a_i$ keeps track of the stabilizer of $\calC/C$ along $\sigma_i$, and we have simple extensions $N_{C, i} \hookrightarrow \pi_*N_{\calC, i}$ corresponding to multiplication by $a_i$ on characteristic monoids. 

\begin{definition} Let $(\calX, M_{\calX})$ be a separated Deligne--Mumford stack with fs log structure. A \emph{log twisted prestable map} is a diagram
\[ \begin{tikzcd}
(\calC, M_\calC) \arrow{r}{(f, f^\flat)}  \arrow{d} & (\calX, M_{\calX}) \\ 
  (S, M_S)\rlap{,}
\end{tikzcd} \]
where $\calC \to S$ is a twisted curve and $(\calC, M_{\calC}) \to (S, M_S)$ is log smooth. A log twisted prestable map is \emph{stable} if the underlying map $f$ is representable and the coarse map $h \colon  C/S \to X$ is a stable map. 
\end{definition}

We will also make use of the \emph{root stack} construction of Borne--Vistoli~\cite[Section 4.2]{bornevistoli} in the case of a locally free log structure. Suppose $(X,M_X)$ is a locally free log scheme and $\phi \colon  \overline{M}_X \to \overline{M}'$ is a simple extension. The root stack $\sqrt[\phi]{(X,M_X)}$ is a log algebraic stack over $\operatorname{Sch}_X$ 
which associates to $s \colon  T \to X$ the groupoid of pairs $(l, \varphi)$, where $l \colon  s^*M_X \to M_T'$ is a simple extension of log structures and $\varphi$ is an isomorphism of maps of characteristic monoids $\overline{l} \cong s^*\phi$. 

We conclude this section by comparing the sheaves of log differentials between $\calC$ and $C$. This simple comparison is one of the advantages of using log geometry to study twisted curves. Let $(f,f^\flat) \colon  (\calC, M_{\calC}) \to (S, M_S')$ and $(h, h^\flat) \colon  (C, M_C) \to (S, M_S)$ denote the natural log smooth maps, where $M_S \hookrightarrow M_S'$ is a simple extension. 

\begin{notat} We set $\Omega^{log}_{\calC/S} := \Omega^{\log}_{(f,f^\flat)}$ and $\Omega^{\log}_{C/S} : = \Omega^{\log}_{(h,h^\flat)}$ and similarly for the cotangent complex.
\end{notat}

\begin{proposition}
Let $\pi\colon  \calC \to C$ be the coarse space of a twisted curve. Then 
\[
\bLl_{\calC/S} \simeq_{\qis} \Omega^{\log}_{\calC/S}[0] \cong \pi^*\Omega^{\log}_{C/S}[0] \simeq_{\qis} L\pi^*\bLl_{C/S}.
\]

\end{proposition} 

\begin{proof} Since both $(f,f^\flat)$ and $(h,h^\flat)$ are log smooth,  $\bLl_{\calC/S} \simeq_{\qis} \Omega^{\log}_{\calC/S}$ and $\bLl_{C/S}\simeq_{\qis} \Omega^{\log}_{C/S}$. Since $\Omega^{\log}_{C/S}$ is locally free, $L\pi^*\Omega^{\log}_{C/S}[0] \simeq_{\qis} \pi^*\Omega^{\log}_{C/S}[0]$. Therefore, it suffices to show that the natural map
$$
\rho \colon  \pi^*\Omega^{\log}_{C/S} \to \Omega^{\log}_{\calC/S}
$$
is an isomorphism which we can compute \'etale locally. Toward that end, suppose $S = \Spec(R)$. In an \'etale neighborhood of a marked gerbe with  stabilizer $\mu_a$, the map $\pi$ can be written as $R[u] \mapsto R[x]$ with $u \mapsto x^a$, so $\rho \colon  d \log u \mapsto a d\log x$. In a neighborhood of a node with stabilizer $\mu_a$, we have $R[u,v]/(uv - t^a) \to R[x,y]/(xy - t)$ with $(u,v) \mapsto (x^a, y^a)$, so $\rho$ is given by $d\log u \mapsto a d\log x$, $d \log v \mapsto a d \log y$. In either case, $\rho$ is an isomorphism since $a \in R$ is invertible. 
\end{proof}

\section{Log maps to 1-dimensional targets}

In this section, we collect some results on log twisted maps to stacky curves which we will use to lift the smoothings from Section~\ref{sec:combinatorics} to the twisted case.

\subsection{Introducing log structures}

We begin by giving a criterion to lift a prestable map with $1$-dimensional target to a log map. More precisely, suppose $g \colon  C \to X$ is a prestable map to a smooth curve $X$ with log structure $M_X$ such that $(X,M_X)$ is log smooth over the trivial log point. We wish to find a log structure $M_C$ on $C$ such that $(C,M_C)$ is a log smooth curve and $g$ lifts to a log map. For the purpose of constructing smoothings, it suffices by Lemma~\ref{lem:comb} to consider the case where the connected components of $C$ which are contracted to the marked points of $(X,M_X)$ are smooth (see also Remark~\ref{rem:f*M}). 

\begin{proposition}\label{prop:liftmap} Let $g \colon  C \to X$ be a prestable map with $(X, M_X)$ as above. Suppose that the $1$-dimensional connected components of the support of $f^{*}(\overline{M}_X)$ are smooth genus $0$ curves. Then $g$ lifts to a  prestable logarithmic map. That is, there exist log structures $M_C$ on $C$ and $M$ on $S = \Spec(k)$ and a map $g^\flat \colon  g^*M_X \to M_C$ such that $(C,M_C) \to (S,M)$ is log smooth. 
\end{proposition}

\begin{proof} The question is \'etale local on the target, so we may pull back to a chart for the log structure $(X,M_X)$. Hence without loss of generality, we may assume that the target is $(\mathbb{A}^1,0) = (X,D)$ with its toric log structure and that $g \colon  C \to \mathbb{A}^1$ is a proper map satisfying the assumptions of the proposition. Moreover, we write $g^{-1}(0) = \{p_1, \ldots, p_k, E_1, \ldots, E_r\}$ as a disjoint union of points and smooth rational curves. Let us denote by $q_{ij}$ the nodes of $C$ along $E_i$ and by $e_{ij}$ the ramification of $g$ at these points. First we claim that $g$ factors through an expansion of the target 
$$
\xymatrix{C \ar[r]^h & \widetilde{X} \ar[r]^\varphi & X.}
$$
Indeed, consider $\widetilde{X} = E \cup_0 \mathbb{A}^1$, where $E = \mathbb{P}(N_{D/X} + \calO_D)$. To construct such a factorization, it suffices to construct non-constant maps $h_i \colon  E_i \to E$ with $h_i(q_{ij}) = 0$ and such that the order of tangency of $h_i$ at $q_{ij}$ is equal to $e_{ij}$. Since the $E_i$ are smooth rational curves, such maps exist. 

Therefore, we have a factorization $h \colon  C \to \widetilde{X}$, and it suffices to show that $h$ can be lifted to a log map where 
$$
M_{\widetilde{X}} = \widetilde{M}_{\widetilde{X}} \oplus_{\calO_{\widetilde{X}}^*} \varphi^*M_X
$$
since $\varphi$ underlies a log map $(\widetilde{X}, M_{\widetilde{X}}) \to (X, M_X)$ and a composition of log maps is a log map. Now we have reduced to the case where the support of $h^*(\overline{M}_{\widetilde{X}}/\mathbb{N})$ is $\{p_1, \ldots, p_k, p_{k + 1}, \ldots, p_n\}$, a union of isolated points. Here the copy of $\mathbb{N}$ in the quotient is the pullback of the minimal log structure on $\Spec(k)$ along $\tilde{X} \to \Spec(k)$. In a neighborhood of the smooth $p_i$, there is a unique divisorial log structure such that the map lifts, encoding the order of tangency. In a neighborhood of the nodal $p_i$, the map lifts to a log map by the tangency condition above, which is precisely the predeformability condition for a map to an expansion (see for example~\cite[Section 12]{HoM}). \end{proof}

\begin{remark} \label{rem:f*M} Proposition~\ref{prop:liftmap} may be generalized to the case where the $1$-dimensional connected components of the support of $f^*\overline{M}_X$ are trees of rational curves at the expense of allowing $C$ to be blown up at the nodes. To do this, one needs to consider the map $g^{\trop}$ of tropical curves and subdivide the source and target to make $g^{\trop}$ tropically transverse (see~\cite[Sections 2.5 and 2.6]{dhruv} for more details). Such a subdivision corresponds to choosing a possibly larger expansion $\widetilde{X} \to X$ and a sequence of log blowups $\widetilde{C} \to C$ after which we may lift $g$ to a map $h \colon  \widetilde{C} \to \widetilde{X}$ which underlies a log map as in the proposition. Since this is not strictly necessary for our stated goal of constructing smoothings, we leave the details to the reader. 
\end{remark}

Next we address the question of lifting the log structure to a twisted map from its coarse space. More precisely, consider a proper morphism $f \colon  \calC/S \to \calX$  from a twisted curve to a smooth Deligne--Mumford curve $\calX$. Let $M_\calX$ be a log structure making $(\calX, M_{\calX})$ a log twisted curve with coarse space $(X, M_X \hookrightarrow M'_X)$. Here $(X,M_X)$ is log smooth, and $M_X \hookrightarrow M'_X$ is the simple extension encoding the stabilizers of $\calX$. Suppose furthermore that there exist log structures $M_C$ and $M_S$ on the coarse space $C$ of $\calC$ and on $S$, respectively, such that
$$
\xymatrix{(C, M_C) \ar[r]^g \ar[d] & (X, M_X) \\ (S, M_S) & }
$$
is a prestable log map. 

\begin{proposition}\label{prop:rootstack} Let $f \colon  \calC/S \to \calX$ and $f \colon  (C,M_C)/(S, M_S) \to (X,M_X)$ be as above. Then there exist a root stack $\pi \colon  \tilde{\calC} \to \calC$ and log structures $M_{\tilde{\calC}}$ on $\tilde{\calC}$ and $M_S'$ on $S$ such that $(\tilde{\calC}, M_{\tilde{\calC}}) \to (S, M_S')$ is a log smooth curve, and there exists a map  $f^\flat \colon  (f \circ \pi)^*M_\calX \to M_{\tilde{\calC}}$ making $f$ into a prestable log map.   
\end{proposition} 

\begin{proof} Let $g \colon  C \to X$ be the coarse map. By Theorem~\ref{thm:olsson}, we have simple extensions of log structures $M_C \hookrightarrow M'_C$, $M_S \hookrightarrow M'_S$, and $M_X \hookrightarrow M_X'$ on $C$, $S$, and $X$, respectively, such that $(C \to S, \{a_i, \sigma_i\}, M_S \hookrightarrow M'_S)$ and $(X \to \Spec(k), \{b_i, x_i\}, k^* \to k^*)$ are log twisted curves. Here the $\sigma_i \colon  S \to C$ correspond to the marked gerbes of $\calC \to S$ with order of stabilizer given by $a_i$. Similarly, the $x_i \in X$ correspond to the marked gerbes of $\calX$ with stabilizer order $b_i$. 

The question is local on $X$, so we can restrict to a neighborhood of one of the marked points $x_i$ with stabilizer $\mu_{b_i}$ on $X$ and suppose without loss of generality that there is a single marked point $x \in X$ with stabilizer $\mu_b$. This point is exactly the support of $\overline{M}_X = \mathbb{N}$ and $\overline{M}_X' = \mathbb{N}$, and the simple extension is simply the map $1 \mapsto b$. 

By assumption, there exists a map $g^\flat \colon  g^*M_X \to M_C$ making $(g, g^\flat)$ into a log map. By Theorem~\ref{thm:olsson}, it suffices to show that there exist a further simple extension $M_C' \hookrightarrow M_C''$ which is an isomorphism away from finitely many smooth points and a map of log structures $f^\flat \colon  g^*M_X' \to M_C'$ such that the square
\vskip -1em
$$
\xymatrix{g^*M_X \ar[r]^{g^\flat} \ar[d] & M_C \ar[d]^{s} \\ g^*M_X' \ar[r]_{f^\flat} & M_C''
}
$$
commutes. Note that $M'_X = M_X \oplus \mathbb{N}/(e = be')$, where $e$ is the element corresponding to a local parameter of $x \in X$ and $e'$ is the generator of the copy of $\mathbb{N}$. By the universal property of this pushout, to construct $f^\flat$, it suffices to show that $s(g^\flat(e))$ is divisible by $b$. However, we can guarantee this by taking a large enough simple extension $s$.  \end{proof} 

\subsection{Deformations of log twisted maps and their coarse spaces}\label{sec:logdef}

Let
\vskip -2em
$$
\xymatrix{(\calC, M_{\calC}) \ar[r]^{(f,f^\flat)} \ar[d] & (\calX, M_{\calX}) \\ (S, M_S)}
$$
be a stable log map from a log twisted curve, and suppose that $(\calX, M_{\calX})$ is a log smooth curve. We define the \emph{critical locus} of $f$ to be the support of the cokernel of the natural map
$$
 f^*\Omega^{\log}_{\calX} \to \Omega^{\mathrm{log}}_{\calC/S}.
$$
Note that this map, viewed as a two-term complex, is a presentation of the log cotangent complex $\bLl_f$ relative to $S$, so we may equivalently define the critical locus to be the support of $\calH^0(\bLl_f)$. 

The critical locus is a union of nodes of $\calC$, branch points of $f$, and contracted components of $f$. By contrast, the kernel of the above map, or equivalently $\calH^{-1}(\bLl_f)$, is supported solely on the contracted components of $f$. 

Let $A_1, \ldots, A_k$ be the connected components of the critical locus and $U_1, \ldots, U_k$ be affine \'etale neighborhoods of the $A_i$ with log structure pulled back from $\calC$. Suppose each $U_i$ avoids $A_j$ for $j \neq i$. Let $\mathrm{Def}_{(f,f^\flat)}$ be the versal deformation space of the log map $(f,f^{\flat})$. Since the moduli space is a Deligne--Mumford stack, $\mathrm{Def}_{(f,f^{\flat})}$ is pro-representable. Let $\mathrm{Def}_i$ be the miniversal deformation space of the restriction $(f_i := f|_{U_i} \colon  U_i \to \calX, f^{\flat}|_{U_i})$. We have the following proposition (see also~\cite[Proposition 4.3]{Vakil}).

\begin{prop}\label{prop:defisom} In the setting above, the  natural map 
$$
\Def_{(f,f^{\flat})} \to \Def_1 \times \cdots \times \Def_k
$$
\vskip -.5em
\noindent is an isomorphism.
\end{prop}

\begin{proof}
The deformation space $\Def_{(f,f^\flat)}$ (resp.\ $\Def_i$) is constructed via the complex 
$$
R\Hom(\bLl_f, \calO_\calC) \ (\text{resp. } R\Hom(\bLl_{f_i}, \calO_{U_i})).
$$

Since the $g_i \colon  U_i \to \calC$ are strict \'etale,  $Lg_i^*\bLl_f = \bLl_{f_i}$. Let $K_i$ be the component of $K = \calH^{-1}(\bLl_f)$ supported on $A_i \subset U_i$ and $Q_i$ the component of $Q = \calH^0(\bLl_f)$ supported on $A_i \subset U_i$. Since the $A_i$ are disjoint, we have $K = \oplus K_i$ and $Q = \oplus Q_i$ and $K_i = \calH^{-1}(\bLl_{f_i})$ and $Q_i = \calH^{0}(\bLl_{f_i})$.

Since both $\bLl_f$ and $\bLl_{f_i}$ are two-term complexes, we have distinguished triangles
$$
K[1] \to \bLl_f \to Q \to\quad\text{and}\quad
K_i[1] \to \bLl_{f_i} \to Q_i \to.
$$

There is a natural map $R\Hom(\bLl_f, \calO_\calC) \to R\Hom(\bLl_{f_i}, \calO_{U_i})$ given by restricting $f$ to $U_i$. This map corresponds to the map on deformation spaces $\Def_{(f,f^\flat)} \to \Def_i$. Taking direct sums yields the  diagram of distinguished triangles
\begin{equation}\label{commutativediagram}
\xymatrix{
R\Hom(Q, \calO_{\calC}) \ar[r] \ar@{=}[d] & R\Hom(\bLl_f, \calO_{\calC}) \ar[r] \ar[d] & R\Hom(K[1], \calO_{\calC}) \ar[r] \ar@{=}[d] & \\
\bigoplus R\Hom(Q_i, \calO_{U_i}) \ar[r]  & \bigoplus R\Hom(\bLl_{f_i}, \calO_{U_i}) \ar[r] & \bigoplus R\Hom(K_i[1], \calO_{U_i})\ar[r] & 
}
\end{equation}
where the vertical maps are equalities via the identifications $K = \oplus K_i$ and $Q = \oplus Q_i$. The middle map corresponds to the product of the restriction map on deformation spaces
$$
\Def_{(f,f^{\flat})} \to \Def_1 \times \cdots \times \Def_k. 
$$
Since the vertical maps on the ends are equalities, the middle map is an isomorphism by the properties of distinguished triangles, or equivalently the five lemma, from which we conclude. \end{proof}

\begin{remark}\label{rmk:criticallocus}
Via Diagram \eqref{commutativediagram}, we see that in fact the deformation spaces $\Def_i$ of an \'etale neighborhood $U_i$ of the component $A_i$ of the critical locus is independent of the choice of $U_i$ since the $Q_i$- and $K_i$-terms are independent of $U_i$. 
\end{remark} 

Now we address the question of lifting log smooth deformations along the coarse space map of a log twisted curve. Let 
$$
\xymatrix{(C_0, M_{C_0}) \ar[r]^{h_0} \ar[d] & (X, M_X) \\ (S_0, M_{S_0}) & }
$$
be a prestable log map, and suppose that there exists a simple extension of log structures $M_S \to M'_{S}$ such that $(C_0/S_0, (\sigma_i, a_i), l_0 \colon  M_{S_0} \to M'_{S_0})$ is a log twisted curve. Let $(\calC_0, \Sigma_i) \to S_0$ be the corresponding twisted curve with log map $(\calC_0, M_{\calC_0}) \to (S_0, M'_{S_0})$, and let $g_0 \colon  (\calC_0, M_{\calC_0}) \to (X, M_X)$ be the composition of $h_0$ with the coarse space map $\pi_0\colon  \calC_0 \to C$. Let $(S_0, l_0 \colon  M_{S_0} \to M'_{S_0}) \hookrightarrow (S, l \colon  M_S \to M'_S)$ be a strict closed immersion of log schemes with simple extensions defined by a square $0$ ideal $I \subset \calO_S$. 

\begin{prop}\label{prop:logflatdef}  With notation as above, suppose we have a log smooth deformation
\begin{equation}\label{eqn:def}
\xymatrix{(C, M_{C}) \ar[r]^{h} \ar[d] & (X, M_X) \\ (S, M_{S}) & }
\end{equation}
of\, $h_0$ over $(S,M_S)$. Then $h$ lifts uniquely to a log smooth deformation 
$$
\xymatrix{(\calC, M_{\calC}) \ar[r]^{g} \ar[d] & (X, M_X) \\ (S, M'_{S}) &  }
$$
of $g_0$ factoring through $h$.
\end{prop} 

\begin{proof} We can summarize the situation with the commutative diagram below:

\[
\xymatrix{
(\calC_0, M_{\calC_0}) \ar[rd]^{\pi_0} \ar@/^/[rrrd]^{g_0} \ar@/_/[ddr] & & & \\
& (C_0, (\pi_0)_*M_{\calC_0}) \ar[r]^{(\mathrm{id}, s_0)} \ar[d] & (C_0, M_{C_0}) \ar[d] \ar[r]_{h_0} & (X, M_X) \\
& (S_0, M'_{S_0}) \ar[r]^{(\mathrm{id}, l_0)} & (S_0, M_{S_0})\rlap{,} & 
}\]
where the middle square is a pullback in the category of fs log schemes by~\cite[Lemma 4.6, Diagram 4.6.2, and Lemma 4.7]{Olsson}. Given a strict square $0$ thickening $(S, l \colon  M_S \to M'_S)$ of $(S_0, l_0 \colon  M_{S_0} \to M'_{S_0})$ with ideal $I$ and a log smooth deformation $h$ as in Diagram~\eqref{eqn:def}, we can form the pullback 
$$
\xymatrix{(C, M'_C) \ar[r] \ar[d] & (C, M_C) \ar[d] \ar[r]^h & (X, M_X) \\ 
(S, M'_S) \ar[r] & (S, M_S) & }
$$
in the category of fs log schemes to obtain a log smooth deformation of $h_0 \circ (\mathrm{id}, s_0)$. Then $M_C \hookrightarrow M'_C$ is a simple extension, and by Theorem~\ref{thm:olsson}, this is equivalent to a log smooth deformation $g$. \end{proof}

\section{Smoothing twisted stable maps}

Our goal now is to lift the smoothings constructed in Section~\ref{sec:combinatorics} to the case of twisted maps. More precisely, suppose we have a genus $0$ twisted stable map $(f_0\colon \calC_0 \to \calX, p_1, \ldots, p_n)$ such that 
\begin{itemize} 
\item $\calX$ is a smooth stacky curve with coarse map $\tau \colon  \calX \to X = \mathbb{P}^1$, 
\item the coarse map $(h_0\colon  C_0 \to X, q_1, \ldots, q_n)$ satisfies Conditions~\ref{conditions} relative to $\infty \in X$, and
\item $\tau$ is \'etale over $\infty$.
\end{itemize}
The situation can be summarized in the following diagram, where $S_0 = \Spec(k)$:
$$
\xymatrix{
(\calC_0, p_1, \ldots, p_n) \ar[r]^-{f_0} \ar[d]_{\pi_0} & \calX \ar[d]^\tau \\ 
(C_0, q_1, \ldots, q_n) \ar[r]^-{h_0} \ar[d] & X \\
S_0\rlap{.} & }
$$

By Theorem~\ref{thm:sufficient}, there is a smoothing $h \colon  (C, \sigma_i)/S \to X$ of $h_0$ over $S = \Spec(R)$, the spectrum of a DVR, with generic fiber contained in $\mapsd$, and we  wish to lift this to a smoothing $f \colon  (\calC, \Sigma_i)/S \to \calX$ of $f_0$ with generic fiber contained in $\calM_{\Gamma}(\calX)$. That is, we wish to fill in the dotted arrows in the diagram   
$$
\xymatrix{
(\calC_0, p_i) \ar@{-->}[r] \ar[d]_{\pi_0} & (\calC, \Sigma_i) \ar@{-->}[r]^-{f} \ar@{-->}[d]_\pi & \calX \ar[d]^\tau \\ 
(C_0, q_i) \ar[r] \ar[d] & (C, \sigma_i) \ar[r]^-{h} \ar[d] & X \\
S_0 \ar[r] & S\rlap{.} & }
$$

The difficulty is that there are global obstructions to lifting from a map to $X$ to a map to $\calX$. For example, the $j$-invariant of a map to $\overline{\calM}_{1,1}$ must satisfy that the discriminant is a sum of a square and a cube which cuts out a high-codimension locus inside the space of maps to $\mathbb{P}^1$. The key observation is that $\tau$ is \'etale. Thus we can lift such a smoothing locally around $\infty$ (see Theorem~\ref{thm:local}) to obtain a partial smoothing of $f_0$. Then we can construct a further global smoothing when $\calX$ has some positivity (see Theorem~\ref{thm:smoothing1}).

\begin{theorem}\label{thm:local} In the situation above, there exists a partial smoothing $f \colon  (\calC, \Sigma_i)/S \to \calX$ of $f_0$ such that 
\begin{itemize} 
\item $f(\Sigma_i) = \infty$,
\item $\calC_\eta$ is smooth in a neighborhood of $f_\eta^{-1}(\infty)$, and 
\item the coarse map on the generic fiber $h_\eta$ is ramified to order $d_i$ at $(\sigma_i)_\eta$. 

\end{itemize}
\end{theorem}

\begin{proof} We wish to use the log deformation theory results of Section~\ref{sec:logdef}. To do this, we begin by showing that the smoothing of the coarse map $h \colon  (C, \sigma_i)/S \to X$ can be endowed with the structure of a stable log map. 

\begin{lemma}\label{log:smoothing} Let $(h_0 \colon  C_0 \to \bP^1, q_1, \ldots, q_n)$ be as above, and let
\[ \begin{tikzcd}
C \arrow{r}{h}  \arrow[d] & \bP^1 \\ 
  S \arrow[u, bend left, "\sigma_i"] 
\end{tikzcd} \] 
be a smoothing constructed in Theorem~\ref{thm:sufficient}. Denote by $M$ the divisorial log structure on $\bP^1$ corresponding to the point $\infty$. Then there exist log structures $M_C$ and $M_S$ on $C$ and $S$ and a map $h^\flat \colon  h^*M \to M_C$ such that $h$ is a stable log map over $S$.
\end{lemma}

\begin{proof} First note that $C \to S$ is log smooth, where we equip $S = \Spec(R)$ with the standard log structure and $C$ with the divisorial log structure induced by the union of the central fiber and marked sections. Indeed, this follows since $C \to S$ is a toroidal morphism of the corresponding toroidal embeddings. Endow $\bP^1$ with the divisorial log structure for $\infty \subset \bP^1$;   the pullback of the monomial $x_\infty$ cutting out $\infty$ is locally a sum of monomials for the toroidal structure on $C$ since $h^*[\infty]$ is supported on the toroidal boundary. Thus, there is a natural map $h^\flat$ of divisorial log structures recording the exponents of the monomial $h^*x_\infty$. \end{proof}

Suppose $(R,\mathfrak{m})$ is complete, and let $h_n\colon  (C_n, \sigma_{i,n})/S_n \to X$ be the compatible system of deformations over $S_n = \Spec (R/\mathfrak{m}^n)$ obtained by truncating the smoothing over $S$. We will lift these to a compatible system of deformations of $h_0$ in two steps.

First let $(g_0, g_0^\flat) \colon  (\calC_0, M_{\calC_0}) \to (X, M)$ be the composition $h_0 \circ \pi_0$, where we note that $h_0$ is a log map by Lemma~\ref{log:smoothing} and $\pi_0$ is a log map by the definition of log twisted curves. By Lemma~\ref{prop:logflatdef} and induction on $n$, the map $h_n$ extends to a log smooth deformation $g_n \colon  (\calC_n, M_{\calC_n}) \to (X, M_X)$ over $(S_n, M'_{S_n})$, where $M_S \hookrightarrow M_S'$ a simple extension and $(S_n, M'_{S_n}) \hookrightarrow (S, M_S')$ is a strict closed embedding. Moreover, the $g_i$ are compatible by construction. In this way, we get an $R$-point of the formal deformation space $\Spf R \to \Def_{(g_0, g_0^\flat)}$.  

Next we will use Proposition~\ref{prop:defisom} to lift this $R$-point of $\Def_{(g_0, g_0^\flat)}$, at least in a neighborhood of $\infty$, to a smoothing of $f_0$. Note that $f_0$ lifts to a log map $(f_0, f_0^\flat)$ by Proposition~\ref{prop:rootstack}, where we endow $\calX$ with the pullback log structure $\tau^*M$. Equivalently, this is the divisorial log structure for $\infty$ but is different from the canonical log structure on a log twisted curve. Let $U \subset \mathbb{P}^1$ be an \'etale neighborhood of $\infty$ such that $\tau|_V \colon  V \to U$ is \'etale, where $V = \tau^{-1}(U)$, and such that $U$ avoids all the components of the critical loci of both $f_0$ and $g_0$ away from $\infty$. Let $\Def_{f_0, \infty}$ be the miniversal deformation space of the restriction of $(f_0, f_0^\flat)$ to an \'etale neighborhood of $f_0^{-1}(\infty)$ and similarly for $\Def_{g_0, \infty}$. Note that by the independence of the chosen \'etale neighborhood (see Remark~\ref{rmk:criticallocus}), we have that these are the miniversal deformation spaces of $(f_0)|_{U} \colon  f_0^{-1}(U) \to U \subset \mathbb{P}^1$ and $(g_0)|_{V} \colon  g_0^{-1}(V) \to V \subset \calX$. We set $U' = f_0^{-1}(U) = g_0^{-1}(V)$ and denote by $U''$ a neighborhood of the critical loci of both $f_0$ and $g_0$ away from $\infty$.

By Proposition~\ref{prop:defisom}, we have
$$
\Def_{(f_0, f_0^\flat)} \cong \Def_{f_0, \infty} \times \Def_{f_0, \neq \infty}
$$
and
$$
\Def_{(g_0, g_0^\flat)} \cong \Def_{g_0, \infty} \times \Def_{g_0, \neq \infty}, 
$$
where $\Def_{f_0, \neq \infty}$ (resp.\ $\Def_{g_0, \neq \infty}$) is the miniversal deformation space of a neighborhood of the critical loci of $(f_0,f_0^\flat)$ (resp.\ $(g_0, g_0^\flat)$) which are away from $\infty$.

\begin{lemma}\label{lem:naturalmap}
The natural map $\Def_{f_0, \infty} \to \Def_{g_0, \infty}$ induced by composition with $\tau$ is an isomorphism.
\end{lemma}

\begin{proof} We proceed as in Proposition~\ref{prop:defisom} and see that there is a natural isomorphism 
$$
R\Hom(\bLl_{f_0}, \calO_{\calC}) = R\Hom(\bLl_{f_0'}, \calO_{U'}) \oplus R\Hom(\bLl_{f_0''}, \calO_{U''})
$$
and similarly
$$
R\Hom(\bLl_{g_0}, \calO_{\calC}) = R\Hom(\bLl_{g_0'}, \calO_{U'}) \oplus R\Hom(\bLl_{g_0''}, \calO_{U''}).
$$

We now compute
$$
\bLl_{f_0'} = [f_0^*\Omega_{U} \to \Omega^{\log}_{U'}],
$$
$$
\bLl_{g_0'} = [g_0^*\Omega_{V} \to \Omega^{\log}_{U'}].
$$

\noindent Since $\tau|_V$ is \'etale, we have that $g_0^*\Omega_V \cong g_0^*\tau^*\Omega_U = f_0^*\Omega_U$ by the functorality of pullbacks. Therefore, the natural map
$$
\bLl_{g_0'} \to \bLl_{f_0'}
$$
is an isomorphism, so we conclude that the natural map
$$
R\Hom(\bLl_{f_0'}, \calO_{U'}) \to R\Hom(\bLl_{g_0'}, \calO_{U'})
$$
is an isomorphism.  \end{proof}

Now let $\Spf R \to \Def_{f_0, \infty}$ be the composition of the formal $R$-point of $\Def_{(g_0, g_0^\flat)}$ constructed above with composition
$$
\Def_{(g_0, g_0^\flat)} \to \Def_{g_0, \infty} \to \Def_{f_0, \infty}, 
$$
where the first map is the projection from Proposition~\ref{prop:defisom} and the second map is the inverse of the isomorphism in Lemma~\ref{lem:naturalmap}. There is a canonical splitting of the projection $\Def_{(f_0, f_0^\flat)} \to \Def_{f_0, \infty}$ via the isomorphism of Proposition~\ref{prop:defisom} given by picking the constant deformation in $\Def_{f_0, \neq \infty}$. Further composing with this splitting gives us a formal deformation $\Spf R \to \Def_{(f_0, f_0^\flat)}$. Since $\Def_{(f_0, f_0^\flat)}$ is the formal neighborhood of a point of a locally of finite type algebraic stack (see~\cite{Chen, GS}), this formal deformation algebraizes to a family $(f,f^\flat) \colon  (\calC, M_{\calC})/(S, M_S') \to (\calX, \tau^*M)$ of log twisted stable maps (see for example~\cite[Theorem 1.1]{Bhatt} and~\cite[Corollary 1.5]{BDHL}). Let us denote by $f \colon  (\calC, \Sigma_i)/S \to \calX$ the underlying twisted stable map. By construction, the coarse map of $f$ agrees with $h$ over an \'etale neighborhood of $\infty$, so we conclude that $f$ satisfies the required conditions. \end{proof}

\begin{remark} \label{rem:away} By construction, the smoothing of Theorem~\ref{thm:local} induces the constant infinitesimal deformation of the restriction of $f_0$ away from $\infty$. More precisely, suppose $A$ is a component of the critical locus of $f_0$ which is disjoint from $f_0^{-1}(\infty)$, and let $U \to C$ be an \'etale neighborhood of $A$ which is disjoint from all other critical loci. Then $f|_{U} \mod \mathfrak{m}^n$ is the constant infinitesimal deformation of $f_0|_{U \times_C C_0}$. 
\end{remark}

\begin{remark}\label{rem:rootstack} In Theorem~\ref{thm:local} and the situation above, we can drop the assumption that $\tau$ is \'etale over $\infty$ at the cost of replacing $\calC_0$ with a root stack. Indeed, if $\tau$ is not \'etale, then $\tau$ factors through a non-trivial root stack along $\infty$. In this case, one obtains a corresponding simple extension of log structures $M \to M'$ on $\bP^1$. Then by Proposition~\ref{prop:rootstack}, the map $f_0$ can be extended to a log map after taking a further root stack $\calC_1 \to \calC_0$. The rest of the argument goes through as written to produce a partial smoothing of the composition $\calC_1 \to \calC_0 \to \calX$ with the required properties as in the theorem. We leave the details to the reader since our main interest is in the case of $(\overline{\calM}_{1,1}, \infty)$, where $\tau$ is \'etale over $\infty$. 
\end{remark}

\subsection{Global smoothings for target weighted projective lines}

Next we consider the problem of smoothing twisted stable maps to a weighted projective line $\calP = \calP(a,b)$ without any tangency conditions.

\begin{theorem}\label{thm:smoothing1} Let $(f_0 \colon  \calC_0 \to \calP, \Sigma_{i,0})$ be a genus $0$ twisted stable map to a weighted projective line $\calP = \calP(a,b)$. Then there exists a smoothing 
\[ \begin{tikzcd}
\calC \arrow{r}{f}  \arrow[d] & \calP \\ 
  S \arrow[u, bend left, "\Sigma_i"] 
\end{tikzcd} \] 
over $S = \Spec (R)$, the spectrum of a DVR.  
\end{theorem} 

\begin{proof} Suppose $(R,\mathfrak{m})$ is complete. We will construct such a smoothing by building a compatible system of deformations over $R/\mathfrak{m}^n$. Since the stack of not necessarily representable twisted stable maps to $\calP$ is a locally of finite type algebraic stack, such a compatible system automatically algebraizes as in the proof of Theorem~\ref{thm:local}. First note that if $a = ka'$ and $b = kb'$ with $a', b'$ coprime, then the canonical map $\rho \colon  \calP(a,b) \to \calP(a',b')$ is \'etale, and so infinitesimal deformations lift uniquely along $\rho$. Thus, without loss of generality, we may assume that $a$ and $b$ are coprime and $\calP$ is a twisted curve. 

  First we reduce to the case where every connected component in the fibers of $f_0 \colon  \calC_0 \to \calP$ is irreducible. We can argue as in Lemma~\ref{lem:comb}. If $E$ is a connected component of the preimage $f_0^{-1}(p)$ for some $p$, then we can view $f_0|_E$ as a map to the residual gerbe $\calG_p$ where $E$ is pointed by the nodes. Then since $\calG_p$ is \'etale over $\Spec(k)$, it suffices to smooth the composition $E \to \Spec(k)$ as a marked twisted curve and then glue together the smoothing of $E$ with the  constant family of maps $f|_{\overline{C_0 \setminus E}}$ along the markings corresponding to nodes exactly as in the proof of Lemma~\ref{lem:comb} to obtain a partial smoothing where each connected component of the contracted locus is irreducible. Therefore, without loss of generality, we may suppose that every connected component of the contracted locus of $f_0$ is irreducible.

By Propositions~\ref{prop:liftmap} and~\ref{prop:rootstack}, there exists a twisted curve $\tilde{\calC}_0$ with a partial coarse space map $\pi_0 \colon  \tilde{\calC}_0 \to \calC_0$ such that $g_0 := f_0 \circ \pi_0 \colon  \tilde{\calC}_0 \to \calP$ can be endowed with a log structure where $\calP$ has the toric log structure and $(\tilde{\calC}_0, M_{\tilde{\calC}_0}) \to (\Spec(k), M)$ is log smooth. Since $(\calP, M_{\calP})$ is toric, the log cotangent complex $\bLl_{\calP} \cong \calO_\calP$ is trivial. Moreover, since $\tilde{\calC}_0$ is genus $0$, 
\begin{equation}\label{eqn:obs}
\Ext^1(g_0^*\bLl_{\calP}, \calO_{\tilde{\calC}_0}) =  H^1(\tilde{\calC}_0, \calO_{\tilde{\calC}_0}) = 0.
\end{equation}
Therefore, the log deformation space of $(g_0,g_0^\flat)$ is log smooth over the log deformation space of the curve $(\tilde{\calC}_0, M_{\tilde{\calC}_0})$. By Theorem~\ref{thm:olsson}, $\Def(\tilde{\calC}_0, M_{\tilde{\calC}_0})$ is a root stack over $\Def(C_0, p_1, \ldots, p_n)$, where the $p_i$ are the smooth marked points lying under $\Sigma_{i,0}$. In particular, we can build a smoothing $(C/S, \sigma_i)$ of $(C_0, p_i)$ whose truncations $(C_n/\Spec (R/\mathfrak{m}^n), \sigma_{i,n})$ lift to a formal smoothing of $(\tilde{\calC}_0,M_{\tilde{\calC}_0})$, after possibly extending the log structure on~$S$. By the vanishing obstruction group (\ref{eqn:obs}), the map $(g_0, g_0^\flat)$ can then be lifted to a compatible family $(g_n \colon  \tilde{C}_n \to \calP, g_n^\flat)$ of log smooth deformations over $R/\mathfrak{m}^n$. 

Thus, after algebraizing, we have a smoothing $(g \colon  \tilde{\calC}/S \to \calP, \Sigma_i)$ of $g_0$. Taking the relative coarse space of~$g$ gives a flat family of maps $f \colon  \calC/\Spec (R) \to \calP$ whose formation commutes with base change (see~\cite[Proposition 3.4]{AOV11} and~\cite[Corollary 3.3]{AOV08}). Therefore, we conclude that the central fiber of $f$ agrees with $f_0$ and that the generic fiber is smooth since it is a partial coarse space of the smooth generic fiber of~$\tilde{C}/S$.   
\end{proof} 

\begin{remark}\label{rem:needlog} For the proof of Theorem~\ref{thm:smoothing1}, it is necessary to use log deformation theory even though the statement does not make any reference to log geometry. The reason is that the stack $\calP$ is \emph{not} convex in the classical sense. The tangent bundle is not necessarily effective even though it has positive degree, and there are examples of genus $0$ twisted stable maps where the vanishing (\ref{eqn:obs}) does not hold for the usual tangent bundle. For example, consider $\calC = [E/\tau]$, where $E$ is an elliptic curve with hyperelliptic involution $\tau \colon  E \to E$, $\pi \colon  \calC \to \bP^1$ is the coarse space, and $p_i \in \bP^1$ are the points lying under the fixed points of $\tau$. Then the line bundle $L = \calO_{\calC}\left(\frac{1}{2}p_1 - \frac{1}{2}p_2 + \frac{1}{2}p_3 - \frac{1}{2}p_4\right)$ with sections $0 \in L$ and $1 \in L^{\otimes 2}$ induces a map $f \colon  \calC \to \calP(1,2)$ with $f^*T_{\calP} = L$ and $h^1(\calC, L) = h^1(\bP^1, \pi_*L) = h^1(\bP^1, \calO(-2)) \neq 0$. 
\end{remark}

\section{Proofs of the main theorems}

We use the notation from the introduction. Let $(\calX, x_1, \ldots, x_r)$ be a smooth and proper pointed stacky curve, and suppose that the coarse map $\tau \colon  \calX \to X$ is \'etale at the points $x_i$ (see also Remark~\ref{rem:rootstack}). Fix a tuple of discrete data $\Gamma = (n_0, \{\Gamma_1, x_1\}, \ldots, \{\Gamma_r, x_r\})$. 

\begin{proposition} The loci $\calM_\Gamma(\calX)$ and $\calN_{\Gamma}(\calX)$ $($see Definition~\ref{def:MGamma}\,$)$ are locally closed substacks of\, $\calK_{0,n}(\calX, d)$. 
\end{proposition} 

\begin{proof} We run through the conditions of Definition~\ref{def:MGamma}. Condition~\eqref{def:MGamma-1} is both open and closed by the definition of a twisted curve, and condition~\eqref{def:MGamma-2} is closed. Conditions~\eqref{def:MGamma-4} and~\eqref{def:MGamma-4'} are open by the openness of the smooth locus. For condition~\eqref{def:MGamma-3}, we consider the cohomology sheaves of the cotangent complex $\bLl_f$. The locus where $f$ is constant is given by the support of $\calH^{-1}(\bLl_f)$. Since the universal family of curves over the moduli space is proper, the locus where the marked points avoid $\Supp \calH^{-1}(\bLl_f)$ is open, and therefore so is the condition of $f$ being non-constant along the marked points.
  Finally, under this condition, the order of tangency of $f$ at $x_i$ is measured by the length of $\mathrm{coker}(f^*\Omega^1_\calX \to \Omega^1_\calC)_x$, and by semicontinuity, there is a locally closed subset where this length is constant $d_i$. 
\end{proof} 

\begin{proof}[Proof of Theorem~\ref{mainthm}]  The marked points $p_{0k}$ for $k = 1, \ldots, n_0$ with no tangency conditions do not change the outcome of the problem. Indeed, if we can construct such a smoothing with $n_0 > 0$, then by forgetting and stabilizing, we obtain a smoothing with $n_0 = 0$, and, conversely, if we have a partial 1-parameter smoothing over a base  $S$ with $n_0  = 0$, then up to a finite base change $S' \to S$, we can pick generic sections passing through $p_{0k}$ to obtain a smoothing with $n_0 > 0$. So without loss of generality, suppose $n_0 = 0$. 

We proceed by induction on $r$. The base case is precisely Theorem~\ref{thm:local}. In general, let $(f \colon  \calC \to \calX, \{\{p_{jk}\}_{k = 1}^{n_j} \}_{j = 1}^r)$ be a map in $\calK_{\Gamma}(\calX)$. Applying Theorem~\ref{thm:local} to $(f \colon  \calC \to \calX, p_{r1}, \ldots, p_{rn_r})$ produces a \mbox{$1$-parameter} deformation whose generic fiber is a relative map to $(\calX, x_r)$ with tangency $\Gamma_r$ which also satisfies that the curve is smooth along the preimage of $x_r$. Moreover, by Remark~\ref{rem:away}, the generic fiber is contained in $\calK_\Gamma(\calX)$. Therefore, without loss of generality, we may assume that $\calC$ is smooth in a neighborhood of $f^{-1}(x_r)$. By the inductive hypothesis, there exists a 1-parameter deformation of $f$ with generic fiber contained in $\calN_{\Gamma'}(\calX)$, where $\Gamma' = (\{\Gamma_1, x_1\}, \ldots, \{\Gamma_{r-1}, x_{r-1}\})$. Moreover, by Remark~\ref{rem:away}, this partial smoothing induces the constant deformation in an \'etale neighborhood of $f^{-1}(x_r)$, and so the generic fiber is contained in~$\calN_{\Gamma}(\calX)$. 
\end{proof} 

\begin{proof}[Proof of Theorem~\ref{mainthm2}] Let $(f \colon  \calC \to \calX, p_{jk})$ be a map contained in $\calK_\Gamma(\calX)$. By Theorem~\ref{mainthm}, there exists a partial smoothing of $(f \colon  \calC \to \calX, p_{jk})$ with generic fiber contained in $\calN_\Gamma(\calX)$, so without loss of generality, suppose that $(f \colon  \calC \to \calX, p_{jk})$ is contained in $\calN_\Gamma(\calX)$. By Theorem~\ref{thm:smoothing1}, there exists a smoothing of the stable map into the interior of $\calK_{0,n}(\calX, d)$ since $\calX = \calP$ is a weighted projective line. Let $\Spf R \to \Def_f$ be the corresponding formal deformation. By Proposition~\ref{prop:defisom}, we can write
$$
\Def_f = \Def_{f, x_1} \times \cdots\times \Def_{f,x_r} \times \Def_{f, \neq}\,,
$$
where $\Def_{f,x_i}$ is a miniversal deformation space of a small \'etale neighborhood of the fiber $f^{-1}(x_i)$ and $\Def_{f, \neq}$ is a miniversal deformation space of the critical loci of $f$ which are disjoint from the fibers $f^{-1}(x_i)$ for each $i$. Then projecting onto $\Def_{f,\neq}$ and then composing with the section $\Def_{f,\neq} \to \Def_f$ which picks the constant deformation on each of the $\Def_{f,x_i}$ factors produces a new formal deformation $\Spf R \to \Def_f$ which agrees with the smoothing from Theorem~\ref{thm:smoothing1} away from the fibers $f^{-1}(x_i)$ but is the constant deformation in a neighborhood of $f^{-1}(x_i)$ for all $i$. After algebraizing this formal deformation as in the proof of Theorem~\ref{thm:local}, we obtain a 1-parameter family over $\Spec (R)$ such that the generic fiber is contained in $\calN_{\Gamma}(\calX)$ but is also smooth away from the union of the fibers $f^{-1}(x_i)$. Therefore, the generic fiber is smooth everywhere and thus contained in $\calM_\Gamma(\calX)$. \end{proof}

\newcommand{\etalchar}[1]{$^{#1}$}
\providecommand{\bysame}{\leavevmode\hbox to3em{\hrulefill}\thinspace}
\providecommand{\MR}{\relax\ifhmode\unskip\space\fi MR }
\providecommand{\MRhref}[2]{%
  \href{http://www.ams.org/mathscinet-getitem?mr=#1}{#2}
}
\providecommand{\href}[2]{#2}

\end{document}